\documentstyle[11pt]{article} 
\begin{document}
\def\r{\rightarrow}
\newcommand{\fdem}{\hfill$\Box$}

\newcommand{\E}{{I\!\!E}}     
\renewcommand{\P}{{I\!\!P}}     
\renewcommand{\L}{{I\!\!L}}     
\newcommand{\F}{{I\!\!F}}     
\newcommand{\T}{\mbox{\rule{0.1ex}{1.5ex}\hspace{0.2ex}\rule{0.1ex}{1.5ex}
\hspace{-1.7ex}\rule[1.5ex]{1.6ex}{0.1ex} }}
     
\newcommand{\N}{{I\!\!N}}     
\newcommand{\Z}{{Z\!\!\!Z}}
\newcommand{\R}{{I\!\!R}}     
\newcommand{\D}{{I\!\!D}}     
\newcommand{\C}{\mbox{\rm C\hspace{-1.1ex}\rule[0.3ex]{0.1ex}{1.2ex}
\hspace{0.3ex}}}
\newcommand{\Q}{\mbox{\rm Q\hspace{-1.1ex}\rule[0.3ex]{0.1ex}{1.2ex}
\hspace{1.1ex}}}

\renewcommand{\Im}{\mbox{\rm Im }}
\renewcommand{\Re}{\mbox{\rm Re }}
\newcommand{\supp}{\mathop{\rm supp}}
\newcommand{\diag}{\mathop{\rm diag}}
\renewcommand{\deg}{\mathop{\rm deg}}
\renewcommand{\dim}{\mathop{\rm dim}}
\renewcommand{\ker}{\mathop{\rm Ker}}
\newcommand{\id}{\mbox{\rm1\hspace{-.2ex}\rule{.1ex}{1.44ex}}
   \hspace{-.82ex}\rule[-.01ex]{1.07ex}{.1ex}\hspace{.2ex}}
\newcommand{\bra}{\langle\,}
\newcommand{\ket}{\,\rangle}
\newcommand{\obl}{/\!/}
\newcommand{\mapdown}[1]{\vbox{\vskip 4.25pt\hbox{\bigg\downarrow
  \rlap{$\vcenter{\hbox{$#1$}}$}}\vskip 1pt}}
\newcommand{\tr}{\mbox{\rm tr}}
\newcommand{\Tr}{\mbox{\rm Tr}}
\newcommand{\r}{\mathop{\rightarrow}}
\newcommand{\cA}{\mbox{$\cal A$}}
\newcommand{\cB}{\mbox{$\cal B$}}
\newcommand{\cC}{\mbox{$\cal C$}}
\newcommand{\cD}{\mbox{$\cal D$}}
\newcommand{\cE}{\mbox{$\cal E$}}
\newcommand{\cF}{\mbox{$\cal F$}}
\newcommand{\cG}{\mbox{$\cal G$}}
\newcommand{\cH}{\mbox{$\cal H$}}
\newcommand{\cJ}{\mbox{$\cal J$}}
\newcommand{\cK}{\mbox{$\cal K$}}
\newcommand{\cL}{\mbox{$\cal L$}}
\newcommand{\cM}{\mbox{$\cal M$}}
\newcommand{\cN}{\mbox{$\cal N$}}
\newcommand{\cP}{\mbox{$\cal P$}}
\newcommand{\cR}{\mbox{$\cal R$}}
\newcommand{\cS}{\mbox{$\cal S$}}
\newcommand{\cT}{\mbox{$\cal T$}}
\newcommand{\cW}{\mbox{$\cal W$}}
\newcommand{\cY}{\mbox{$\cal Y$}}
\newcommand{\cX}{\mbox{$\cal X$}}
\newcommand{\cU}{\mbox{$\cal U$}}
\newcommand{\cV}{\mbox{$\cal V$}}

\centerline{\bf STABLE LAWS}
\medskip
\centerline{\bf  AND PRODUCTS OF POSITIVE RANDOM MATRICES}
\bigskip
\centerline{\bf --------------------------------}
\bigskip
\centerline{{\Large\bf H. Hennion} \footnote{Hubert Hennion, \ Institut Math\'ematiques de Rennes, 
Universit\'e de Rennes I, Campus de Beaulieu, 35042 Rennes-Cedex, France, 
\ Email : Hubert.Hennion@univ-rennes1.fr.}, 
{\Large\bf L. Herv\'e} \footnote{Loic Herv\'e, I.R.M.A.R, Institut National des Sciences Appliqu\'ees, Campus de Beaulieu, 
35042 Rennes-Cedex, France,\ Email : Loic.Herve@insa-rennes.fr}} \medskip
\centerline{\bf------}
\vskip 1truecm
\noindent
{\bf Summary} \par
    Let $S$ be the multiplicative semigroup of $q\times q$ matrices 
with positive entries such that every row and every column contains 
a strictly positive element. 
Denote by $\displaystyle (X_n)_{n\geq1}$ a sequence of independent identically 
distributed random variables in $S$ and by
$\displaystyle X^{(n)} = X_n \cdots X_1$, $\, n\geq 1$,
the associated left random walk on $S$. 
We assume that $\displaystyle (X_n)_{n\geq1}$ verifies the contraction property
\par\noindent
\centerline{$\displaystyle \P\Bigl(\bigcup_{n\geq1}[X^{(n)} \in S^\circ]\Bigr)>0$,}
where $S^\circ $ is the subset of all matrices which have strictly positive
entries. 
We state conditions on the distribution of the random matrix $X_1$ 
%which are similar to the ones used in  the case of real random variables,
which ensure that the logarithms of the entries, of the norm, 
and of the spectral radius of the products $X^{(n)}$, $n\ge 1$,
are in the domain of attraction of  a stable law.
%We prove that, if the logarithm of the entries of  $X_1$ are in the 
%domain of attraction of a stable law, then  
%the sequences of the logarithm of the entries 
%and of the spectral radius of 
%the random products $\displaystyle X^{(n)} = X_n \cdots X_1$, $\, n\geq 1$,
%when suitably normalized, converge to this stable law .\par  
\vfill\eject
\noindent{\bf I. STATEMENT OF THE RESULT}

\noindent Let $S$ be the multiplicative semigroup of $q\times q$ matrices 
with real non negative entries such that every row and every column 
contains a strictly positive element. The subset of $S$ composed of 
matrices with strictly positive entries is a subsemigroup of $S$
denoted by $S^\circ$. \\
\noindent Let $(e_i)_{i=1,\ldots ,q}$ be the canonical basis of the linear 
space $\R^q$. Then a $q\times q$ matrix is identified with 
an endomorphism of $\R^q$. 
We denote by $\langle \cdot, \cdot \rangle $ the canonical  scalar product  on $\R^q$, and 
we define the cones $C$ and $\overline{C}$ by
$$C=\{x : x\in \R^q, \forall i=1,\ldots ,q, \, \langle  x,e_i \rangle >0 \},
\ \  
\overline{C}=\{x : x\in \R^q, \forall i=1,\ldots ,q, \, \langle  x, e_i \rangle  \geq 0 \}.$$
If $g$ is a $q\times q$ matrix, $g^*$ will stand for its adjoint. We have $g\in S$ [resp. $g\in S^\circ$]
if and only if $g(C)\subset C$ and $g^*(C)\subset C$ [resp. either $g(\overline{C}\backslash\{0\})\subset C$ or 
$g^*(\overline{C}\backslash\{0\})\subset C$]. \\[0.15cm]
The product of $g$ and $g'$ in $S$ is denoted by $gg'$, and for $x\in \overline{C}$, 
$gx$ is the image of $x$ under $g$. Finally $\R^q$ is endowed with the norm $\|\cdot\|$ defined by \\
\centerline{$\displaystyle x\in\R^q,\ \ \ \|x\|=\sum_{i=1}^q|\langle x,e_i\rangle|$.}\\[0.1cm]
Let $(X_n)_{n\geq1}$ be a sequence of independent identically 
distributed (i.i.d) random variables (r.v) in $S$ defined on a probability space 
$(\Omega, \cF ,\P)$. We consider the left random walk 
$$\displaystyle X^{(n)},\ n\geq 1,\ \ X^{(1)}(\omega)=X_1(\omega), \ \ \ X^{(n+1)}(\omega)= X_{n+1} (\omega) X^{(n)} (\omega).$$
Our basic assumption is that $(X_n)_{n\geq1}$ verifies the {\bf contraction property} \\[0.28cm]
$\displaystyle (\cC) \ \ \ \ \ \ \ \ \ \ \ \ \ \ \ \ \ \ \ \ \ \ \ \ \ \ \ \ \ \ \ \ \ \ \ \ \ \ \ \ \ \  \ \ \ 
\P\Bigl(\bigcup_{n\geq 1} [X^{(n)}\in S^\circ ]\Bigr)>0$.\\[0.25cm]
The subsemigroup  $ S^\circ$ is in fact an ideal of $S$, 
that is : if $g\in  S^\circ$ and $g'\in S$,
then $g'g$ and $gg'\in S^\circ$. Consequently $ S^\circ$ is stochastically closed
for the random walk $(X^{(n)})_{n\geq1}$. We set
$$T(\omega)=\inf \{n : n\geq 1, \, X^{(n)}(\omega)\in  S^\circ \}.$$
It is easily shown, Lemma II.1, that~: $(\cC)\ \Leftrightarrow\ \P[T<+\infty ]=1\ \Leftrightarrow\  
\P(\cup_{n\geq 1} [X^{(n)}\in S^\circ ]) = 1$. \\[0.15cm] 
Our aim is to present conditions on $X_1$ ensuring the distributional convergence to a stable law for   
the sequences of real random variables
$$\bigl(1_{[T\le n]}\ln \langle  y, X^{(n)} x \rangle \bigr)_{n\geq 1}, 
\ \ x, y\in \overline{C}\backslash \{0\}.$$
Denoting by $\vec{\bf 1}$ the vector in $\R^q$ whose all entries equal $1$,
we point out that the scalar products $\langle  y, X^{(n)} x \rangle $,
$x, y\in \overline{C}\backslash \{0\}$, include : \\[0.12cm] 
\textindent{-} the matrix entries~:\
$\langle  e_i , X^{(n)} e_j \rangle , \ \ i, j=1,\ldots ,q$, \\[0.12cm]
\textindent{-} the  norm of the image under $X^{(n)}$ of any $x\in \overline{C}\backslash \{0\}$~:
\ $\|X^{(n)} x\|= \langle  \vec{\bf 1}, X^{(n)} x \rangle $, \\[0.12cm]
\textindent{-} the norm $|||X^{(n)} |||=  \langle \vec{\bf 1} , X^{(n)} \vec{\bf 1} \rangle $
of $X^{(n)}$. \\[0.2cm]
Closely related to these quantities is
the spectral radius $\Lambda_n$ of the matrix $X^{(n)}$.
Actually the Perron-Frobenius Theorem yields $\Lambda_n>0$, and  we shall see that the above mentioned distributional 
convergences also concern the sequence $\bigl(\ln \Lambda_n\bigr)_{n\ge 1}$. 

\noindent To state our result, one needs the two following real r.v~: \\[0.1cm]
\centerline{$\displaystyle N_1= |||X_1||| = \sum_{i,j=1}^q \langle e_i, X_1e_j\rangle,\ \ $ and 
$\displaystyle \ \ V_1 = \min_{i=1,\ldots ,q}\sum_{j=1}^q  \langle e_i, X_1e_j\rangle.$} \\[0.1cm]
$N_1$ takes in account the size of the matrix $X_1$ while $V_1$ measures
the smallness of its lines. 

\noindent{\bf Theorem I.} {\it Assume that $(\cC)$ holds and that there exist a real number $\alpha$, 
$0<\alpha\leq2$, a slowly varying function $L$ which is unbounded in case $\alpha=2$, and finally 
some positive constants $c_+$ and $c_-$ with $c_++c_->0$ such that \\[0.1cm]
\textindent{(i)} $\displaystyle \lim_{u\rightarrow +\infty}{u^\alpha\over L(u)}\, \P[N_1> e^u]=c_+$, 
$\displaystyle\ \ \ \lim_{u\rightarrow +\infty}{u^\alpha\over L(u)}\, \P[N_1\le e^{-u}] = c_-$, \\[0.1cm]
\textindent{(ii)} $\displaystyle \limsup_{u\rightarrow +\infty}{u^\alpha\over L(u)}\, \P[V_1\le e^{-u}]<+\infty$.\\[0.2cm]
Then there exist a sequence $(a_n)_{n\ge 1}$
in $\R_+^*$ with $\lim_n a_n=+\infty$ and a sequence 
$(b_n)_{n\ge 1}$ in $\R$ such that,
for any sequences $(x_n)_{n\geq 1}$ and
$(y_n)_{n\geq 1}$ of unit vectors of $\overline{C}$, the random sequences
$$\biggl({1\over a_n}\Bigl( 1_{[T\leq n]}\ln \langle  y_n , X^{(n)} x_n \rangle 
-b_n\Bigr)  \biggr)_{n\geq 1}\ \ \ \ \mbox{and}\ \ \ \ \ 
\biggl({1\over a_n} \Bigl( \ln \Lambda_n -b_n \Bigr) \biggr)
_{n\geq 1}$$
converge in distribution to a stable law of index $\alpha$.} 

\noindent Observe that Hypothesis {\it (i)} means that the real r.v. $\ln N_1$ belongs to 
the domain of attraction of a stable law of index $\alpha$, $0<\alpha\leq2$, the standard Gaussian case being here excluded since 
$L$ is assumed to be unbounded when $\alpha=2$. As it will be seen later on, the hypotheses of Theorem I imply that the above 
considered sequences of random variables have the same distributional behaviour that a sum of i.i.d random variables (See $\S$ IV). 
However it is worth noticing that this is not true when $\alpha=2$ and $L$ is bounded. In fact, to complete the Gaussian case, 
recall it is proved in \cite{hen} that, if Conditions {\it (i)} and {\it (ii)} are replaced by the moment condition 
$\E[|\ln N_1|^2]+\E[|\ln V_1|^2]<+\infty$, then the random sequences of Theorem I converge to a normal law. 
The method used in \cite{hen} is based on martingale techniques, and 
the central limit theorem proved there is also valid when $(X_n)_n$ is supposed stationary and 
satisfies suitable mixing conditions. By the way, recall that, in some cases, the unnormalized random products $(X^{(n)})_{n\geq1}$ 
converge in distribution, see \cite{hen} \cite{kesspi} \cite{muk}. \\[0.2cm] 
Consider the case $q=1$. Then Theorem I corresponds to the well-known statement of convergence to 
stable laws for i.i.d random variables~: we have $S=\R_+^*$, Condition $(\cC)$ holds, and Condition {\it (i)} 
states that $\ln X_1$ is in the above described domain of attraction, {\it (ii)} is a consequence of {\it (i)}. 
So Theorem I gives the expected conclusion. \\[0.2cm] 
%\noindent As briefly explained below, the conclusion in Theorem I will stem from the distributional convergence of a certain 
%i.i.d sequence of real r.v belonging to the above mentioned domain of attraction. Consequently, 
%the sequences $(a_n)_{n\ge 1}$, $(b_n)_{n\ge 1}$, and the limit characteristic function in Theorem I, 
%may be chosen according to the usual theory of stable laws, see \cite{ibra} completed by \cite{aar1}. \\[0.2cm] 
%in particular $(a_n)_{n\ge 1}$ is such that ${n\over a_n^\alpha}L(a_n)=1$, and 
%the limit characteristic function $\Phi$ in Theorem I is defined as follows. by using the theory of stable laws \cite{ibra} 
%completed by \cite{aar1}. \\[0.2cm] 
%Let $J_\alpha=\int_0^{+\infty} {e^{iu}-\chi(u)\over u^\alpha}du$, 
%where, if $0<\alpha<1$, $\chi=0$, \ if $\alpha=1$,
%$\chi=1_{[0,1]}$, and, if $1<\alpha<2$, $\chi=1_{\R_+}$. (It is known that, for $0<\alpha<1$, $J_\alpha=\Gamma(1-\alpha)
%e^{i(1-\alpha)\pi/2}$, while $\Im J_1={\pi\over 2}$.) Then, by defining on $\R^*$ the function 
%$\displaystyle c_\alpha(t) = (c_++c_-)\, \mbox{Im}(J_\alpha)-i{t\over |t|}(c_+-c_-)\, \mbox{Re}(J_\alpha)$, one gets~: \\
%$\displaystyle\Phi(t)=\exp\Bigl(-c_\alpha(t)\, |t|^\alpha\Bigr)$ for $\alpha\not=1\ \ \ \ \ \ ; \ \ \ \ \ \ 
%\displaystyle \Phi(t)= \exp\Bigl(-c_1(t)\, |t|-i(c_+-c_-)\, t \ln|t|\Bigr)$ for $\alpha=1$. 
\noindent The proof of Theorem I is based on the spectral method that was introduced by 
Nagaev \cite{nag1}, \cite{nag2} and later developped by several authors, 
see \cite{hulo}. Although this method has been essentially used to prove Central Limit Theorems and their refinements, we mention that 
Nagaev himself \cite{nag1} has considered the convergence to stable laws, and that his method 
has been extended to the context of dynamical 
systems, See e.g \cite{gui} \cite{aar2} \cite{aar3} \cite{bapei} \cite{gou}. \\[0.2cm]
\noindent Section II summarizes some statements of \cite{hen}, based on Condition $(\cC)$ and related to the projective action 
defined by $g\cdot x={gx \over \|gx\|}$ for $g\in S$ and unit vector $x$ in $\overline{C}$.   
In Section III, denoting by $Y_k$ the adjoint matrix of $X_k$, and setting $\xi (g,x)=\ln \|gx\|$, 
we show that the distributional convergences of Theorem I 
are valid if, for any  unit vector $y$ in $\overline{C}$, the same 
holds for the random variables $\xi(Y_k,(Y_{k-1}\ldots Y_1)\cdot  y)$. Since these r.v may be seen as a functional of the Markov 
chain $(Y_k,(Y_{k-1}\ldots Y_1)\cdot y)_k$, Nagaev's method applies. In fact, it shall be applied to the transition probability 
$P$ of the simpler Markov chain $(Y_k\ldots Y_1\cdot y)_k$, and we shall prove, by using the contractivity properties 
stated in Section II, 
that $P$ satisfies a strong ergodicity condition on a 
certain Lipschitz function space, and finally, by applying the perturbation theory, that the Fourier kernels $P_t$  
associated to $P$ and $\xi$ inherit near $t=0$ the spectral properties of $P$. \\ 
As usual in Nagaev's method, the previous preparation will show that the desired distributional convergence is based on 
the behaviour of 
the power of the dominating eigenvalue $\lambda(t)$ of $P_t$. Actually, one of the main arguments is Proposition III.1 
which links $\lambda(t)$ with the characteristic function of the r.v $\xi$ under the stationary distribution of 
$(Y_k,(Y_{k-1}\ldots Y_1)\cdot y)_k$. So everything shall turn out as in the i.i.d case, provided that $\xi$ 
belongs to the already mentioned domain of attraction. We shall see in  Section IV that this requirement 
holds under Conditions {\it (i)} and  {\it (ii)}. \\[0.2cm] 
\noindent The above relation between the dominating perturbed eigenvalue of the Fourier kernels 
and the characteristic function of the functional 
under invariant distribution has been already exploited in \cite{gou} \cite{bapei}, and mentioned under a different form in \cite{hulo} 
(Lem. IV.4'). It is worth noticing that such a relation holds whenever the spectral method applies, and that 
it greatly makes easier the use of Nagaev's method when dealing with stable laws excluding the standard Gaussian case ; 
for instance it yields a significant simplification of some proofs in \cite{aar2} \cite{aar3}. 

\noindent{\bf II. CONTRACTIVITY}

\noindent {\bf II.1. Preliminaries.} We set
$$B=C\cap \{x : \, x\in\R ^q, \|x\| =1 \}, \ \ \hbox{and} \ \
\overline{B}=\overline{C}\cap \{x : \, x\in\R ^q, \|x\| =1 \},$$
and we define the adjoint random walk $(Y^{(n)})_{n\geq1}$ of $(X^{(n)})_{n\geq1}$
by
$$Y_n = X_n ^*, \ \ 
\ \ \ \ \ Y^{(n)}={X^{(n)}} ^*=Y_1\cdots Y_n, \ \ \ \ n\ge 1.$$
{\bf Lemma II.1.} {\it $(\cC)$ is equivalent to $\P[T<+\infty]=1$. Let $\omega$ be such that $T(\omega)<+\infty$. Then\\[0.15cm]
\textindent{(i)} for  $n\geq T(\omega)$, \ $X^{(n)}(\omega)\in S^\circ$, \\[0.15cm]
\textindent{(ii)} setting, for $n\ge 1$, $\displaystyle D_n(\omega)=\sup \Bigl\{\,
\bigg| 1_{[T\leq n]}(\omega)\, \ln \langle  y , X^{(n)}(\omega) x \rangle  
-\ln \|Y^{(n)} (\omega) y\| \bigg|\ :\  x, y\in \overline{B} \, \Bigr\},$
we have $\displaystyle \sup_{n\ge 1} D_n(\omega)<+\infty$, \\[0.15cm]
\textindent{(iii)} setting $\chi = \frac{1}{q}\vec{\bf 1}$ and, for $n\ge 1$, $\displaystyle \widetilde  D_n(\omega)=
\bigg| \ln \Lambda_n(\omega) -\ln \|Y^{(n)} (\omega)\chi\| \bigg|,$
we have $\displaystyle \sup_{n\ge 1} \widetilde  D_n(\omega)<+\infty$.}

\noindent{\bf Proof.} Suppose $(\cC)$ holds. Then there exists $k\in\N^*$ such that $p=\P[X^{(k)}\in S^\circ]>0$.
The r.v. \par\noindent
$$T'=\inf\{n : n\geq 1, X_{nk}\ldots  X_{(n-1)k+1}\in S^\circ\}$$
has a geometric distribution with parameter $p$. 
Since $S^\circ$ is an ideal, we have $T\leq kT'$, hence $\P[T<+\infty]=1$. The converse implication is obvious. \\ 
Now let any fixed $\omega\in\Omega$ be such that $T(\omega)<+\infty$. Assertion (i) follows from the fact $S^\circ$ is an ideal. 
To prove (ii), it suffices to establish that $\displaystyle\sup_{n\ge T(\omega)}D_n(\omega)<+\infty$. In the following inequalities, 
one considers any fixed integer $n$ such that $n\geq T(\omega)$. For convenience, $\omega$ will be omitted in most of the 
next computations. 
Let $a, b$ be two strictly positive real numbers such 
that, for $i,j=1,\ldots ,q$, we have $\displaystyle a\leq \langle  e_i , X^{(T)}e_j \rangle \leq b$, 
and let $x$ and $y$ be any elements of $\overline{B}$. Using $\|x\|=\|y\|=1$, we obtain for $i=1,\ldots , q$ \par\noindent
$$a \leq \langle  e_i , X^{(T)} x \rangle  \leq b.$$
%Set $\vec{\bf 1}=(1/q,\ldots ,1/q)\in B$. Let $y\in \overline{B}$, 
(That is, $a\vec{\bf 1} \leq  X^{(T)} x \leq b\vec{\bf 1}$ for the coordinatewise order relation on $\R^q$.)  
Moreover, using the formula $\displaystyle \langle  y, X^{(n)}x \rangle = \langle Y_{T+1}\cdots Y_ny , X^{(T)}x \rangle$, 
one gets successively  
%$\displaystyle a\, \langle Y_{T+1}\cdots Y_ny , \vec{\bf 1} \rangle  \leq \langle  y , X^{(n)}x \rangle   \leq
%b\, \langle  Y_{T+1}\cdots Y_ny , \vec{\bf 1} \rangle$, 
$$a\, \langle Y_{T+1}\dots Y_ny , \vec{\bf 1} \rangle  \leq \langle  y , X^{(n)}x \rangle   \leq
b\, \langle  Y_{T+1}\cdots Y_ny , \vec{\bf 1} \rangle,$$
$$|\,\ln \langle  y , X^{(n)}x \rangle  - \ln \|Y_{T+1}\cdots Y_ny\|\,|
\leq \max \{ |\ln a|, |\ln b| \}.$$
In particular, with $x=\vec{\bf 1}$, this gives $|\,\ln\|Y^{(n)}y\| -
\ln \|Y_{T+1}\cdots Y_ny\|\,| \leq \max \{ |\ln a|, |\ln b| \}$. These two inequalities imply  
$\displaystyle\sup_{n\ge T(\omega)}D_n(\omega)<+\infty$. \\
To prove (iii), again consider $\omega\in\Omega$ such that $T(\omega)<+\infty$, and 
recall that, from the Perron-Frobenius Theorem, there exists $R_n(\omega)\in\overline{B}$ such that 
$X^{(n)}(\omega)R_n(\omega)=\Lambda_n(\omega) R_n(\omega)$. 
With $x=R_n$ and $y=\chi = \frac{1}{q}\vec{\bf 1}$, Assertion (ii) yields 
$$\bigg|1_{[T\leq n]}\ln  \frac{1}{q} \langle \vec{\bf 1} , X^{(n)}R_n \rangle -\ln \|Y^{(n)}\chi\|\bigg|\le D_n.$$
From $\langle \vec{\bf 1} , X^{(n)}R_n \rangle = \Lambda_n\|R_n\| = \Lambda_n$, it follows that 
$$\widetilde  D_n(\omega) \le D_n(\omega) + 1_{[T> n]}(\omega)\, |\ln\Lambda_n(\omega)|\, + 1_{[T\le n]}(\omega) \ln q.$$
This proves assertion {\it (iii)}.\fdem

\noindent We deduce from the above lemma that, for any sequence $(a_n)_{n\ge 1}$ in $\R_+^*$
such that $\lim_n a_n=+\infty$, we have \ $\displaystyle \lim_n{1\over a_n}D_n=0$ \ and 
\  $\displaystyle \lim_n{1\over a_n}\widetilde  D_n=0$ a.s. \\[0.15cm]
{\it Consequently, the conclusion in Theorem I for 
$(1_{[T\leq n]}\ln \langle  y_n , X^{(n)} x_n \rangle)_n$ and $(\ln \Lambda_n)_n$ will hold if the same is valid for 
$(\ln \|Y^{(n)} y_n\|)_n$ for any sequence $(y_n)_{n\geq 1}$ of vectors of $\overline{B}$. }
%So our problem is reduced to the behaviour in distribution of the 
%sequences $\bigl(\ln \|Y^{(n)} y\|\bigr)_{n\ge 1}$, $y\in \overline B$.

\noindent{\bf II.2. Projective action of positive matrices.} It is well known that the projective action of matrices 
plays a key part in the study of the asymptotic behaviour of 
random invertible matrix products, cf. \cite{bou} for example. 
As shown in \cite{hen}, this is also true in the case of positive 
matrices. About this action, we now recall the facts that we 
shall use throughout, refering to \cite{hen} for more details and for the proofs. \\[0.1cm]
\noindent Consider the subset $\widetilde  {\overline C}$ of the $q$-dimensional 
projective space associated with the cone $\overline C$. 
In other words, $\widetilde {\overline C}$ is the set of lines 
through $0$ and some point in $\overline C\backslash\{0\}$. 
These may be represented by points of the closed polygon $\overline B$. 
An element $g\in S$ maps a line in $\widetilde {\overline C}$ onto 
a line in $\widetilde {\overline C}$, and this defines its projective action on  
$\widetilde {\overline C}$. As $\widetilde {\overline C}$ is represented by 
$\overline B$, the projective action of $g$  moves to the action on 
$\overline B$ defined by  
$$g\cdot x={gx \over \|gx\|}.$$
(Recall that $gx$ is the image of $x$ under the linear action of $g$.) 
The projective action has the following basic properties : if $e$ stands for 
the identity matrix and \ $g, g'\in S$, $x\in \overline B$, we have \par\noindent
\centerline{$e\cdot x=x$, \ \ \ \ $ (gg')\cdot x=g\cdot (g'\cdot x)$.}
%Notice the importance of $\cdot$ to distinguish between the action of 
%elements of $S$ on $\overline C$ and on~$\overline B$. 

\noindent It is well known \cite{bap} that, when $B$ is equipped with the Hilbert 
distance $d_H$, the elements of $S$ have a contractive action 
on $B$,  and that this contractive action is strict for elements
of  $ S^\circ$. However, because the  Hilbert distance is unbounded and only defined on $B$, 
it is more convenient for our purposes 
to use a bounded distance $d$ on $\overline B$
which have similar properties. This distance, already used in \cite{hen}, is defined as follows. 
For $x=(x_1,\ldots ,x_q)$ and $y=(y_1,\ldots ,y_q)$ in $\overline B$, 
we write 
$$m(x,y)=\sup\{\lambda : \lambda\in \R_+, \,\forall i=1,\ldots ,q, 
\, \lambda y_i\leq x_i \}=
\min\{y_i^{-1}x_i : i=1,\ldots ,q, \, y_i>0\}.$$
Besides let $\varphi$ be the one-to-one function on $[0,1]$ defined by $\displaystyle \varphi(s)={1-s \over 1+s}$. 
Then, if $x,y\in\overline B$, one has $\sum_{i=1}^q x_i=\sum_{i=1}^q y_i=1$, thus $0\leq m(x,y)\leq 1$, so one may define 
$$d(x,y)=\varphi\bigl(m(x,y)m(y,x)\bigr).$$
\noindent {\bf Proposition II.1.} (cf. \cite{hen}, $\S$ 10)   
{\it The map $d$ defines a distance on $\overline B$ having the following properties \\[0.1cm]
{\bf (i)} $\sup \{ d(x, y) : x, y \in \overline B \}=1\ \ \ \ $ {\bf (ii)}  if $x, y \in \overline B$, $\|x-y\|\leq 2 d(x,y)$\\[0.1cm]
{\bf (iii)} the topology of $(B, d)$ is the topology induced 
on $B$ by the standart topology of $\R^q$.\\[0.15cm]
Moreover, for $g\in S$, there exists $c(g)$ such that \\[0.1cm]
{\bf (iv)} if $x, y\in \overline B$, $d(g\cdot x, g\cdot y) \leq c(g) d(x, y)\leq c(g)\ \ \ $ 
{\bf (v)} $c(g)\leq 1$, and $c(g)<1$ if and only if $g\in S^\circ$ \\[0.1cm]
{\bf (vi)}  if $g'\in S$, $c(g g')\leq c(g) c(g')\ \ \ $ {\bf (vii)}  $c( g^*)=c(g)$.}

\noindent For any $x\in \overline{B}\backslash B$ and any $y\in B$,
we have $m(x,y)=0$, so that $d(x,y)=1$. Thus \\[0.2cm]
\centerline{$\overline{B}\backslash B = \cup_{x\in \overline{B}\backslash B} \{ y : y\in \overline{B}, d(x,y)<1/2\}$.} \\[0.2cm]
is an open subset of $(\overline B, d)$.
 It follows that the topology of $(\overline B, d)$
and the topology induced by $\R^q$ on $\overline B$  do not coincide ; 
from {\it (ii)} the former is finer than the latter. In the sequel, unless 
otherwise stated, when we appeal to topological properties of 
$\overline B$ and $B$,  we shall assume that these sets are endowed 
with the topologies induced by $\R^q$ ; the distance $d$ will be only used to express contractivity. 

\noindent{\bf II.3. Stochastic contractivity.} Denote by $\mu$ the probability distribution of $Y_1=X_1 ^*$
and by $\mu^{(n)}$ the distribution of $Y^{(n)}=Y_1\ldots Y_n$, 
$n\geq 1$. 
For $n\geq 1$, we set
$$c(\mu^{(n)}) =\sup \Bigl\{\int_S {d(g\cdot y,g\cdot y')\over d(y,y')}
d\mu^{(n)}(g) : y,y'\in \overline B, \, y\not= y' \Bigr\}.$$
Since $c(\cdot)\leq1$, we have $c(\mu^{(n)})\leq1$. Furthermore, the sequence $\displaystyle \bigl(c(\mu^{(n)})\bigr)_{n\geq 1}$ 
is clearly submultiplicative, so we can define
$$\kappa(\mu)=\lim_n c(\mu^{(n)})^{{1\over n}}
=\inf_{n\geq 1} c(\mu^{(n)})^{{1\over n}}.$$
Using Assertion (v) in Proposition II.1, it is easily shown that $(\cC)$ is equivalent 
to $\kappa(\mu)<1$.\par
\medskip\noindent
{\bf Theorem II.1.} {\it Under Condition $(\cC)$, there exists a r.v. $Z_1$ taking values in $B$ 
such that $(Y^{(n)}\cdot x)_n$ converges a.s to $Z_1$, the convergence being 
uniform for $x\in \overline B$. The probability distribution $\nu$ of $Z_1$ verifies $\nu(B)=1$.
It is the unique $\mu$-invariant probability distribution on 
$\overline B$, i.e. the unique probability distribution
on $\overline B$ such that, for any bounded continuous function 
$f$ on $\overline B$, we have }
$$ \int_{\overline B} \Bigl(\int_Gf(g\cdot x)d\mu(g)\Bigr)d\nu(x)
= \int_{\overline B} f(x)d\nu(x).$$
\noindent {\bf Proof.} Using the contractivity properties of $c(\cdot)$,
we see that the sequence of positive r.v. 
$\bigl(c(Y^{(n)})\bigr)_{n\ge 1}$
decreases and hence converges almost surely. Under  $(\cC)$,
there exists an integer $b\in\N^*$ such that $\E[c(Y^{(b)})]<1$. The independence then yields 
$\limsup_k \E[c(Y^{(kb)})]\le\lim_k \bigl(\E[c(Y^{(b)}]\bigr)^k=0$.
It follows from these two facts that 
$\lim_n c(Y^{(n)})=0$ \  a.s.
Notice that, by means of the subadditive ergodic theorem,
we can get, more precisely,  $\displaystyle\lim_n \bigl(c(Y^{(n)})\bigl)^{1\over n}
=\kappa$ \  a.s. \\[0.15cm]
Set $\Omega_1=\{ \omega :\lim_n c(Y^{(n)} (\omega))=0 \}$.
Let $\omega\in\Omega_1$. For $n\ge T(\omega)$, the polygons 
$K_n(\omega)=Y^{(n)} (\omega)\cdot (\overline{B})$ \ form a decreasing 
sequence of compact subsets of $B$, so that
$\displaystyle K(\omega)=\cap_{n\geq 1} K_n(\omega)\not= \emptyset$.
Moreover, for the distance $d$, the diameter $\Delta(\omega)$ of $K(\omega)$ 
is equal to $0$. Indeed we have for $n\geq T(\omega)$,
$$\Delta(\omega)\leq \Delta_n(\omega)=\sup\bigl\{d(Y^{(n)}(\omega)\cdot x , Y^{(n)}(\omega)\cdot y) : 
x, y \in \overline{B} \bigr\}\leq c(Y^{(n)}(\omega)).$$
Define $Z_1(\omega)$ by setting $K(\omega)=\{Z_1(\omega) \}$. Then $Z_1\in K_n(\omega)$ implies 
$d(Y^{(n)}(\omega)\cdot x ,Z_1(\omega)) \leq c(Y^{(n)}(\omega))$, 
and (ii) in Proposition II.1 yields the desired convergence. \\[0.15cm]
%\footnote{????? Compactness then implies that, for any $\varepsilon>0$, there exists a $n_0 (\omega)$ such that, for 
%any $n\ge n_0(\omega)$, $K_n(\omega)\subset \bigl\{y : y\in B, \|y-Z_n(\omega)\|\le \varepsilon\bigr\}.$ } 
\noindent Now denote by $\nu$ the law of $Z_1$.
Set $Z_2=\lim_n (Y_2\ldots Y_n)\cdot x$ \ a.s. 
Clearly $Z_2$ has the distribution $\nu$, and we have 
$Y_1\cdot Z_2=Z_1$ a.s. This gives the $\mu$-invariance
of $\nu$.  Let $\nu'$ be any $\mu$-invariant distribution on 
$\overline B$. Then, for any continuous bounded function 
$f$ on $\overline B$ and $n\ge 1$, we have $\displaystyle \int_{\overline B} \E[f(Y^{(n)}\cdot x)]d\nu'(x)
= \int_{\overline B} f(x) d\nu'(x)$. Thus %$\int_{\overline B} f(x) d\nu(x)=
$\E[f(Z_1)]=\int_{\overline B} f(x) d\nu'(x)$.
Hence $\nu=\nu'$. \fdem

\noindent {\bf III. FOURIER KERNELS }

\noindent{\bf III.1. Definition and link with our distributional problem.} 
Recall that our aim is to study the distributional behaviour of the 
sequences $\bigl(\ln \|Y_1\cdots Y_n y\|\bigr)_{n\ge 1}$, $y\in \overline B$ (cf. the end of $\S$ II.1).
However, since $\bigl((Y_n\ldots  Y_1 \,y)\bigr)_{n\ge 1}$ is
a Markov chain and $Y_n\ldots  Y_1 y$ has the same distribution as $Y_1\cdots Y_ny$, 
it is more convenient to consider
$\bigl(\ln \|Y_n\ldots  Y_1 y\|\bigr)_{n\ge 1}$.
So we introduce  the new left random walk on $S$,
$$\widetilde Y^{(n)}=Y_n\ldots  Y_1,   \ \ n\geq 1,\ \ \widetilde Y^{(0)}=e.$$
For $y_0\in \overline{B}$, consider the sequence of r.v. in 
$\overline{B}$ defined by 
$\bigl(\widetilde  Y^{(n)}\cdot y_0\bigr)_{n\ge 0}$.
It is easily checked that it is a Markov chain on $\overline{B}$ starting at $y_0$
and associated with the transition probability $P$ defined by
$$Pf(x)=\int_S f(g\cdot x)d\mu(g),$$ 
where $x\in \overline{B}$ and $f$ is a bounded measurable function 
on $\overline{B}$. Theorem II.1 shows that $\nu$ is the unique 
$P$-invariant distribution. Finally, for $g\in S$ and $x\in\overline{B}$, define 
$$\xi (g, x)=\ln \|gx\|.$$
The function $\xi$ is connected with the projective action of $S$ on $\overline B$ by the
additive cocycle property
$$\xi (g g',x)=\xi (g,g'\cdot x)+\xi (g',x)\ \ \ (g, g'\in S,\, x\in \overline B).$$
%It will follow, see relation $(\star)$ in Section III below, that the variables 
%under study may be canonically written as sums. 
This property shows that, 
for any $y\in \overline{B}$ and $n\ge 1$, we have 
$$(\star)\ \ \ \ \  \ln\|\widetilde Y^{(n)}y\| =  \xi(\widetilde  Y^{(n)}, y)=
\sum_{k=1}^n \xi(Y_k, \widetilde  Y^{(k-1)}\cdot  y).$$
With the function $\xi$ and the transition probability $P$, we associate the Fourier kernels 
$P_t$, $t\in \R$, 
$$x\in\overline{B},\ \  P_tf(x)=\int_S e^{it\xi(g,x)}f(g\cdot x)d\mu(g),$$
with $f$ as above. The Markov property implies that for $n\ge 1$, $y\in \overline{B}$ and $t\in\R$ (see e.g \cite{hulo}) 
$$(\star\star)\ \ \ \ \ \ \E[e^{it \ln\|\widetilde Y^{(n)}y\|}]={P_t}^n{\bf 1}(y),\ \ \mbox{where}\ \ 
{\bf 1} = 1_{\overline{B}}.$$
This basic relation shows that limit theorems for the 
sequence $\bigl(\xi(\widetilde Y^{(n)}, x)\bigr)_{n\ge1}$ may be 
deduced from the asymptotic behaviour of the iterates of 
the operators $P_t$ acting on a suitable Banach space.
This is the main idea of the spectral method. In Sections III.2-4 below, 
we shall prove that $P$ satisfies a strong ergodicity property on the 
usual space of Lipschitz functions on $\overline{B}$, and we shall apply the standard operator perturbation theorem to 
the Fourier kernels. 

\noindent {\bf III.2. A strong ergodicity property for $P$.} 
We denote by $\cL$ the space of all complex-valued functions $f$ on $\overline{B}$ such that  
$$m(f)=
\sup \Bigl\{{|f(x)-f(x')|\over d(x,x')} : x, x'\in \overline{B},\,  x\not
= x'  \Bigr\}<+\infty.$$
Since the distance $d$ is bounded, the elements of $\cL$
are bounded, so we can equip  $\cL$ with the norm 
$$f\in \cL, \ \ \  \|f\|_{_{\tiny \cL}} = \|f\|_u + m(f), \ \ \
\hbox{with} \ \ \|f\|_u = \sup\{ |f(x)| : x\in \overline{B}\}.$$
Then $(\cL, \|\cdot\|_{_{\tiny \cL}})$ is a  Banach space.
Notice that the functions in $\cL$ may be discontinuous on $\overline{B}$ w.r.t the induced topology of $\R^q$, see the 
remark following Proposition II.1. We still denote by $\|\cdot\|_{_{\tiny \cL}}$ the operator norm on $\cL$, and 
$\Pi$ stands for the rank one projection on $\cL$ defined by~: $\Pi f = \nu(f){\bf 1}$.

\noindent {\bf Theorem III.1.} {\it Under Condition $(\cC)$, for any $\kappa_0\in]\kappa(\mu),1[$, there exists $C>0$ such that, 
for all $n\geq1$, we have $\|P^n-\Pi\|_{_{\tiny \cL}} \leq C\kappa_0^n$. }

\noindent {\bf Proof.} We follow \cite{hulo}. For $x, x'\in M, x\not=x'$, we have  
$${|P^n f(x)-P^n f(x')|\over d(x,x')}\ \leq\ \int{|f(g\cdot x)-f(g\cdot x')|\over d(g\cdot x,g\cdot x')}
{d(g\cdot x,g\cdot x')\over d(x,x')}d\mu^{(n)}(g)\ \leq\  m(f)\,c(\mu^{(n)}),$$
so $P^nf\in\cL$ and $m(P^nf)\le m(f)\,c(\mu^{(n)})$. Since $\|Pf\|_u \leq \|f\|_u$, $P$ acts  continuously on $\cL$. \\
Now set $H=\ker(\nu)\cap\cL$. Since $\nu$ is $P$-invariant and defines a continuous  linear functional on $\cL$, 
$H$ is a closed $P$-invariant subspace in $\cL$. Moreover, when restricted to $H$, the semi-norm $m$ is equivalent to 
the norm $\|\cdot\|_{_{\tiny \cL}}$~: more precisely, if $h\in H$, we have
$$m(h)\leq \|h\|_{_{\tiny \cL}}\leq \bigl(2\sup\{d(y,y') : y,y'\in \overline{B}\}+1\bigr)m(h) \leq 3m(h),$$
the second inequality being deduced from the fact that, if $\nu(f)=0$, there exist $x_1, x_2\in \overline{B}$ such that 
$\Re f(x_1)=\Im f(x_2)=0$. Let $f\in \cL$. Since $f-\Pi(f)\in H$, we have $P^n(f-\Pi(f))\in H$ for all $n\geq1$. Hence 
$$\|P^n(f-\Pi(f))\|_{_{\tiny \cL}} \leq 3\, m(P^n(f-\Pi(f))) \leq 3\, c(\mu^{(n)}) m(f-\Pi(f)) = 
3\, c(\mu^{(n)}) m(f) \leq 3\, c(\mu^{(n)})\|f\|_{_{\tiny \cL}}.$$
Finally, under $(\cC)$, we have $\lim_n c(\mu^{(n)})^{{1\over n}} = \kappa(\mu)$ ($\S$ II.3). 
This gives the desired statement. \fdem 

\noindent {\bf III.3. The Fourier kernels near $0$.} 
To apply the perturbation theory near $t=0$ to the Fourier kernels $P_t$, we have to show that  
$P_t$ is a bounded operator of $\cL$, and to study $\|P_t-P\|_{_{\tiny \cL}}$ when $t\r0$. 
For that, we shall need the following notations. 
For $g\in S$, define 
$$\|g\|=\sup\{\|gx\| : x\in \overline B\}, \ \ v(g)=\inf\{\|gx\| : x\in \overline B\},\ \ \mbox{and}\ 
\ell(g)=|\ln\|g\|\,|+|\ln v(g)|.$$
(Notice that $v(g)>0$.) Finally set $\varepsilon(t)=\int_S \min\bigl\{|t|\ell(g), \, 2\bigr\}d\mu(g)$, 
and observe that $\lim_{t\rightarrow 0}\varepsilon(t)=0$. 

\noindent{\bf Theorem III.2.} {\it For $t\in \R$, $P_t$ defines a bounded operator of $\cL$, and 
$\displaystyle \|P_t-P\|_{_{\tiny \cL}}=O\bigl(\varepsilon(t)+|t|\bigr)$.}
%\centerline{$\displaystyle \lim_{t\rightarrow 0}\|P_t-P\|_{_{\tiny \cL}}=0$.}}

\noindent{\bf Proof.} Recall $P_t$ is associated to $P$ and $\xi(g,x)=\ln \|gx\|$ ($g\in S$, $x\in \overline{B}$). 

\noindent{\bf Lemma III.1.} {\it  For $g\in S$ and $z, x, y\in \overline{B}$ 
such that $d(x,y)<1$, we have }
$$|\xi(g,z)|\le\ell(g), \ \ \ \
| \xi(g, x) -\xi(g,y) |\leq 2 \ln{1\over 1- d(x,y)}.$$
\noindent{\bf Proof.} The first inequality is obvious. The second one is Assertion {\it (ii)} 
of Lemma 5.3 in \cite{hen}, for completeness
we reproduce the proof here. Let $x=(x_1,\ldots ,x_q), \ y=(y_1,\ldots ,y_q)\in \overline{B}$ and 
$g=[g_{ij}]_{i,j=1,\ldots ,q}$. Then 
$\displaystyle \|gx\|=\sum_{i=1}^q \sum_{j=1}^q g_{ij} x_j \geq m(x, y) \sum_{i=1}^q \sum_{j=1}^q g_{ij} y_j=m(x, y)\|gy\|$. 
As $d(x,y)<1$, the number $m(x, y)$ and $m(y,x)$ are in $]0, 1]$. Consequently, the symmetry in $x$ and $y$ yields 
$\displaystyle m(x,y)\leq {\|gx\|  \over  \|gy\| }\leq {1  \over  m(y, x) }$, and 
\begin{eqnarray*}
|\xi(g, x) -\xi(g, y) | \leq \max \{-\ln m(y,x), -\ln m(x,y)\}
&\leq& -\ln m(y,x)-\ln m(x,y)\\
&=& -\ln \varphi^{-1}\bigl(d(x,y)\bigr)=\ln { 1+d(x,y) \over 1-d(x,y)}.
\end{eqnarray*}
For $t\in[0, 1[$, \ $\displaystyle 2\ln {1\over 1-t}-\ln {1+t\over 1-t}
=\ln {1\over 1-t^2} \geq 0$, thus {\it (ii)} follows.\fdem

\noindent Set $\Delta_t=P_t-P$. For $f\in\cL$ and $x\in \overline{B}$, we have 
$\displaystyle \Delta_t f(x) =\int(e^{it\xi(g,x)}-1)f(g\cdot x)d\mu(g)$. \\
Before we proceed, notice the inequality~: $\displaystyle \forall u, v\in\R, \ \ |e^{iu}-e^{iv}|\le \min\{|u-v|, 2\}$.
Thus $$|\Delta_t f(x)|\le \int \min\bigl\{|t|\ell(g), \, 2\bigr\} \, |f(g\cdot x)|d\mu(g)\le \varepsilon(t)\|f\|_u.$$ 
So \ $\|\Delta_t f\|_u\le \varepsilon(t)\, \|f\|_u$. Now for $x, y\in \overline{B}$, 
write $\displaystyle{\Delta_t f(x)-\Delta_tf(y)\over d(x,y)} = A(x,y)+B(x,y)$, with  
$$A(x,y)=\int{e^{it\xi(g,x)}-e^{it\xi(g,y)}\over d(x,y)} f(g\cdot x)d\mu(g)\ \ \mbox{and}\ \ 
B(x,y)=\int \bigl(e^{it\xi(g,y)}-1\bigr){f(g\cdot x)-f(g\cdot y)\over  d(x,y)}d\mu(g).$$
If $d(x,y)>1/2$, we have 
$$|e^{it\xi(g,x)}-e^{it\xi(g,y)}|
\le \min\bigl\{|t| |\xi(g,x)-\xi(g,y)|,2\bigr\}\, (2d(x,y))
\le 4\min\bigl\{|t|\, \ell(g),1\bigr\}\, d(x,y),$$
while,  for $d(x,y)\le 1/2$, the inequality of Lemma III.1 gives
$$|e^{it\xi(g,x)}-e^{it\xi(g,y)}|
\le 2 |t|\ln {1\over 1-d(x,y)}\le 2 C |t| d(x,y),$$
with $C=\sup\{{1\over u}\ln {1\over 1-u} : 0<u \le 1/2\}<+\infty$.
From that, we obtain $\displaystyle |A(x,y)|\le \bigl(4\varepsilon(t)+2C|t|\bigr)\, \|f\|_u$. 
Otherwise, since $c(g)\le 1$, 
$$|B(x,y)|\le \int |e^{it\xi(g,x)}-1|\,
|{f(g\cdot x)-f(g\cdot y)\over  d(g\cdot x,g\cdot y)}|
\ {d(g\cdot x, g\cdot y)\over  d(x,y)}\, d\mu(g)
\le m(f)\,\varepsilon(t).$$
So $\displaystyle m(\Delta_t f)\le \bigl(4\varepsilon(t)+2C|t|\bigr)\,\|f\|_u + \varepsilon(t) \, m(f)$, 
therefore $\|\Delta_t\|_{_{\tiny \cL}}\le 4\varepsilon(t)+2C|t|$.\fdem

\noindent {\bf III.4. Spectral properties of $P_t$ near $t=0$.} 
The following perturbation theorem extends the spectral conclusion of Theorem III.1 to $P_t$ for $t$ near 0.  
Let $\kappa_0$ be chosen as in Theorem III.1. 

\noindent{\bf Theorem III.3.} {\it  We assume that Condition $(\cC)$ holds. Let $\kappa\in]\kappa_0,1[$. 
There exists an open interval $I$ centered at $t=0$ such that, for $t\in I$, 
$P_t$ admits a dominating eigenvalue $\lambda(t)\in\C$, with a corresponding rank-one eigenprojection $\Pi(t)$, satisfying the 
following properties~: }
$$\lim_{t\r0}\lambda(t) = 1,\ \ \|\Pi(t) - \Pi\|_{_{\tiny \cL}} = O(\|P_t-P\|_{_{\tiny \cL}})\ \ \mbox{and}\ \ 
\sup_{t\in I}\|P_t^n -  \lambda(t)^n\, \Pi(t)\|_{_{\tiny \cL}} = O(\kappa^n).$$
\noindent {\it Proof.} We only sketch the proof, refering to \cite{ds} for the details and using standard notations. 
It follows from Theorem III.1 that the spectrum $\sigma(P)$ of $P$ is contained in $\{1\}\cup\overline{D(0,\kappa_0)}$. 
Since $t\mapsto P_t$ is continuous (Th. III.2), there exists $t_0>0$ such that, for $|t|\leq t_0$, we have 
$\sigma(P_t)\subset D(1,\frac{1-\kappa}{2})\cup D(0,\kappa)$, and $\sigma(P_t)\cap D(1,\frac{1-\kappa}{2}) = \{\lambda(t)\}$, where 
$\lambda(t)$ is a simple eigenvalue of $P_t$ with a corresponding rank-one eigenprojectioon $\Pi(t)$ depending continuously on $t$. 
Let $\Gamma$ be the oriented circle $\cC(0,\kappa)$. Since $(z,t)\mapsto (z-P_t)^{-1}$ is continuous on the compact 
set $\Gamma\times[-t_0,t_0]$, the formula $P_t^n-\lambda(t)^n\Pi(t) = \frac{1}{2i\pi}\int_\Gamma z^n (z-P_t)^{-1}dz$ leads to 
the last estimate of Theorem. \fdem

\noindent The next proposition states a simple expansion for the  perturbed eigenvalue $\lambda(t)$. 

\noindent {\bf Proposition III.1.} {\it For $t\in I$, we have $\lambda(t) = \mu\otimes\nu(e^{it\xi}) + O(\|P_t-P\|_{_{\tiny \cL}}^2)$. }

\noindent{\it Proof.} Since $\nu$ defines a continuous  linear functional on $\cL$ and $\|P_t-P\|_{_{\tiny \cL}}\r0$ when $t\r0$, 
the rank-one eigenprojection $\Pi(t)$, defined in Theorem III.3, 
is such that $\nu(\Pi(t){\bf 1}) \r \nu(\Pi{\bf 1}) = 1$. So one may assume that $\nu(\Pi(t){\bf 1})\neq0$ 
for any $t\in I$, with $I$ possibly reduced. For $t\in I$, set $v(t) = (\nu(\Pi(t){\bf 1}))^{-1}\, \Pi(t){\bf 1}$. 
Then we have $\lambda(t)v(t) = P_tv(t)$ and $\nu(v(t)) = 1$, therefore 
$$\lambda(t) = \nu(P_tv(t)) = \nu(P_t {\bf 1}) + \nu(P_t(v(t)-{\bf 1})) = 
\mu\otimes\nu(e^{it\xi}) + \nu((P_t-P)(v(t)-{\bf 1})),$$
the last equality following from $\nu(P(v(t)-{\bf 1})) = \nu(v(t)-{\bf 1}) = 0$ since $\nu$ is $P$-invariant.  
We conclude by observing that $\|v(t)-{\bf 1}\|_{_{\tiny \cL}} = \|v(t)-v(0)\|_{_{\tiny \cL}} = 
O(\|\Pi(t) - \Pi\|_{_{\tiny \cL}}) =  O(\|P_t-P\|_{_{\tiny \cL}})$. \fdem

\noindent{\bf IV. PROOF OF THEOREM I}

\noindent Let us point out how the results of the previous sections will be used to establish Theorem I. 
We have to study the distributional behaviour of \\
\indent $\displaystyle \ \ \ \ \ \ \ \ \ \ \ \ \ \ \ \ \ \ \ \ \ \ \ \ \ \ \ \ \ \ \ \ 
\ln\|\widetilde Y^{(n)}y_n\| = \sum_{k=1}^n \xi(Y_k, \widetilde  Y^{(k-1)}\cdot  y_n)$ \\
for any sequence $(y_n)_{n\geq 1}$ of vectors of $\overline{B}$.
Let $y\in\overline{B}$, from Theorem II.1 and the independence of $Y_k$ and 
$\widetilde  Y^{(k-1)}$, the sequence 
$(\xi(Y_k, \widetilde  Y^{(k-1)}\cdot  y))_k$ converges 
in distribution to $\xi(Y_1,Z_2)$, where $Z_2$ is independent of 
$Y_1$ and has distribution $\nu$. From this we may guess that 
$\bigl(\ln \|\widetilde Y^{(n)}y_n\|\bigr)_{n\ge 0}$ has the same asymptotical 
behaviour that a sequence of sums of stationary random variables distributed 
as $\xi(Y_1,Z_2)$. This is confirmed by a look at the characteristic functions.
In fact, the   characteristic function of $\ln \|\widetilde Y^{(n)}y_n\|$
is $P(t)^n1(y_n)$ whose asymptotic behaviour is, as shown by 
the spectral decomposition of Theorem III.3, essentially ruled by 
$\lambda(t)^n$ which doesn't depend on $y_n$ or on any 
initial distribution. Notice that this can be used as will be done in the 
sequel to deduce a limit theorem from the expansion of $\lambda(t)$
at $0$, but also conversely to get an expansion  of $\lambda(t)$
at  $0$ from a known limit theorem, see \cite{ihp1} \cite{ihp2}.  
Now observe that the  characteristic function of $\xi(Y_1,Z_2)$ is 
$\mu\otimes\nu(e^{it\xi})$, which is precisely the first term in the 
expansion of $\lambda(t)$ in Proposition III.1. So we see that if, for a sequence $(a_n)_n$ of positive real numbers, we have  
$$\ (\star \star \star)\ \ \ \ \|P_{\frac{t}{a_n}}-P\|_{_{\tiny \cL}}^2 = o(\frac{1}{n}),$$ 
then $(\mu\otimes\nu(e^{i\frac{t}{a_n}\xi}))^n$ is the principal part of the expansion of $(\lambda(\frac{t}{a_n}))^n$, so that 
$\bigl(\frac{1}{a_n}\ln \|\widetilde Y^{(n)}y_n\|\bigr)_{n\ge 0}$
will behave as $\bigl(\frac{1}{a_n}\sum_{k=1}^n \Xi_k\bigr)_{n\ge 1}$,
where $\bigl(\Xi_k\bigr)_{k\ge 1}$ is a sequence of independent 
random variables distributed as $\xi(Y_1,Z_2)$. \\ 
Actually we shall show that, under the hypotheses of 
Theorem I, the law of $\xi(Y_1,Z_2)$ is in the domain of attraction of a 
stable law and that $(\star \star \star)$ is verified with the corresponding scaling sequence $(a_n)_n$,  
these two facts lead to the claimed result. Finally observe that Condition 
$(\star \star \star)$ is not fulfilled with $a_n=\sqrt n$ in the standard gaussian case since 
it is known that in this case the variance of the limit law is not 
$\sigma^2\bigl( \xi(Y_1,Z_2)\bigr)$ \cite{hulo}. 

%\noindent Untill now we only used the contractivity condition $(\cC)$.  
%Here we shall also need Conditions (i)-(ii) in Theorem I, which concern the r.v 
%$\displaystyle N_1=\sum_{i,j=1}^q <e_i, X_1e_j>$ and $\displaystyle V_1=\min_{i=1,\ldots ,q}\sum_{j=1}^q <e_i, X_1e_j>$. \\
%\noindent For convenience let us rewrite here these two conditions~:  
%{\it there exist $\alpha$, $0<\alpha\leq2$, a slowly varying function $L$ which is unbounded in case $\alpha=2$, and 
%finally some positive constants $c_+$ and $c_-$ with $c_++c_->0$ such that we have, when $u\r+\infty$ \\[0.15cm]
%(i) $\displaystyle {u^\alpha\over L(u)}\, \P[N_1> e^u] \r c_+, \ \ \ {u^\alpha\over L(u)}\, \P[N_1\le e^{-u}] \r c_-\ \ \ ; $ 
%(ii) $\ \ \ \limsup{u^\alpha\over L(u)}\, \P[V_1\le e^{-u}]<+\infty$. } 

\noindent {\bf Proposition IV.1.} {\it  Suppose that Conditions $(\cC)$ and (i)-(ii) hold. Let $F$ be 
the distribution function of $\xi(Y_1,Z_2)$~: 
$$F(u)=\mu\otimes \nu\{ (g,x) : \xi(g,x)\leq u\}\ \ \  (u\in \R).$$
Then there exist positive functions $\rho_+$ and $\rho_-$ 
defined on $\R_+^*$ such that, for $u>0$,}
$$1-F(u)={\rho_+(u)L(u)\over u^\alpha},\ \ \ F(-u)={\rho_-(u)L(u)\over u^\alpha},\ \ \mbox{with}\ \ 
\lim_{u\rightarrow +\infty} \rho_+(u)=c_+\ \ \mbox{and} \ \ \lim_{u\rightarrow +\infty} \rho_-(u)=c_-.$$

\noindent {\bf Proof.} Since $\xi(g,x)=\ln\|gx\|$ and  $\mu$ is the law of $Y_1$, 
we have $\displaystyle F(u) =\int_{\overline{B}}\P[\|Y_1x\|\leq e^u]d\nu(x)$. For $x\in B$, we set 
$\displaystyle N_1^x=\|Y_1x\|=\sum_{i,j=1}^q <e_i, X_1 e_j>x_i$. 
As $\nu(B)=1$, one gets for $u>0$ 
$$1-F(u)=\int_B \P[N_1^x > e^u] d\nu(x)\ \ \ \mbox{and}\ \ \ F(-u)=\int_B \P[N_1^x \leq e^{-u}] d\nu(x) .$$
To proceed, we have to compare the tails of $N_1^x$, for $x\in B$,
with that of $N_1$. Let $u_0$ be such that, for $u\ge u_0$, we have $L(u)>0$. 
For $u>u_0$, and for $x\in B$, we set  \par\noindent
\centerline{$\displaystyle \ \ c_+(u)={u^\alpha\over L(u)}\P[N_1>e^u]$, \ \ 
$\displaystyle c_-(u)={u^\alpha\over L(u)}\P[N_1\le e^{-u}]$,} 
\centerline{$\displaystyle\ \  c_+(x,u)={u^\alpha\over L(u)}\P[N_1^x>e^u]$, \ \ 
$\displaystyle c_-(x,u)={u^\alpha\over L(u)}\P[N_1^x\le e^{-u}]$.} \\[0.2cm]
\noindent  Setting $\displaystyle m(x)=\min_{i=1,\ldots , q} x_i$, we have the following obvious inequalities 
$m(x) N_1\le N_1^x\le N_1$, and 
$$\P[m(x)N_1>e^u]\le \P[N_1^x>e^u]\le  \P[N_1>e^u]$$ 
$$\P[N_1\le e^{-u}]\le \P[N_1^x\le e^{-u}]\le  \P[m(x)N_1\le e^{-u}].$$
Now let $\varepsilon$, $0<\varepsilon<1$. Suppose that $u>0$ 
is such that $e^{-\varepsilon u}\le m(x)$,
we get\par\noindent
\centerline{$\displaystyle
{u^\alpha\over L(u)}c_+\bigl((1+\varepsilon)u\bigr)
{L\bigl((1+\varepsilon)u\bigr)\over \bigl((1+\varepsilon)u\bigr)^\alpha}
\le c_+(x,u)\le c_+(u)$,}  
\centerline{$\displaystyle c_-(u)\le c_-(x,u) \le 
{u^\alpha\over L(u)}c_-\bigl((1-\varepsilon)u\bigr)
{L\bigl((1-\varepsilon)u\bigr)\over \bigl((1-\varepsilon)u\bigr)^\alpha}.$}\\[0.15cm]
Since by hypothesis, $c_+(v)\r c_+$ and $c_-(v)\r c_-$ when $v\r+\infty$, it follows that \par\noindent
\centerline{$\displaystyle {c_+\over (1+\varepsilon)^\alpha}
\le \liminf_{u\rightarrow +\infty}c_+(x,u)
\le \limsup_{u\rightarrow +\infty}c_+(x,u)\le c_+$,}
\centerline{$\displaystyle 
c_-\le \liminf_{u\rightarrow +\infty}c_-(x,u)
\le \limsup_{u\rightarrow +\infty}c_-(x,u)
\le {c_-\over (1-\varepsilon)^\alpha}.$} \\[0.15cm]
Thus $\displaystyle \lim_{u\rightarrow +\infty}c_+(x,u)= c_+ $ and $\displaystyle\lim_{u\rightarrow +\infty}c_-(x,u)= c_-$. \\[0.15cm]
Lastly $N_1^x\le N_1$ yields $\P[N_1^x>e^u]\le  \P[N_1>e^u]$, while 
$N_1^x=\|Y_1x\| = \sum_{j=1}^q \|Y_1e_j\|x_j \geq V_1$ gives $\P[N_1^x\le e^{-u}]\le  \P[V_1\le e^{-u}]$. 
Therefore, for any $x\in B$ and $u>0$, we have 
$$c_+(x,u)\leq c_+(u)\ \ \mbox{and}\ \ c_-(x,u)\leq \frac{u^\alpha}{L(u)} \P[V_1\le e^{-u}],$$ 
and by (i)-(ii), the functions of the variable $u$ in each right term of these inequalities are bounded on $\R_+$. 
Now one may conclude. We have 
$${u^\alpha\over L(u)}\bigl(1-F(u)\bigr) = \int_B c_+(x,u)d\nu(x)\ \ \ \mbox{and}\ \ \ \ 
{u^\alpha\over L(u)}F(-u) = \int_B c_-(x,u)d\nu(x),$$
and Lebesgue's Theorem implies that these integrals converge to $c_+$ and $c_-$ respectively as $u\r+\infty$.~\fdem

\noindent Proposition IV.1 means that $\xi(Y_1,Z_2)$ belongs to the domain of attraction of a stable law with order 
$\alpha$, $0<\alpha\leq2$,  the standard Gaussian case being excluded since, for $\alpha=2$, $L$ is supposed to be unbounded. \\
Let $(\Xi_k)_k$ be an independent sequence of real r.v distributed as 
$\xi(Y_1,Z_2)$. From Proposition IV.I, there exist sequences $(a_n)_n$ in $\R_+^*$ and $(b_n)_n$ in $\R$ such that 
$(\frac{\Xi_1+\ldots+\Xi_n\, -\, b_n}{a_n})_n$ converges in distribution to a stable law of order 
$\alpha$, see \cite{ibra}. It is known that $\lim_na_n=+\infty$ and that 
$(a_n)_{n\ge 1}$ may be chosen such that 
$\displaystyle {n\over a_n^\alpha}L(a_n)=1$. 

\noindent {\bf Proposition IV.2.} {\it  Suppose that Conditions $(\cC)$ and (i)-(ii) hold. For any fixed real $t$, we have 
$\|P_{\frac{t}{a_n}}-P\|_{_{\tiny \cL}}^2 = o(\frac{1}{n})$ when $n\r+\infty$. }

\noindent{\bf Proof.} The notations $\|g\|$, $v(g)$ and $\ell(g)$ below have been introduced for Theorem III.2.  
We have $V_1=v(Y_1)$~: 
indeed, by definition, $V_1=\min_{i=1,\ldots ,q}\|Y_1e_i\|$, so $e_i\in\overline{B}$ implies $v(Y_1) \leq V_1$, and if 
$x\in\overline{B}$, one has $\|Y_1x\| = \sum_{j=1}^q \|Y_1e_j\|x_j \geq V_1$, 
thus $v(Y_1) \geq V_1$. 
Besides $|||g||| = \langle \vec{\bf 1} , g \vec{\bf 1} \rangle$ is a norm for $q\times q$-matrices,  
while, for $g\in S$, the quantity $\|g\|$ corresponds to the matrix norm associated to the norm $\|\cdot\|$ on $\R^q$. 
Since the two previous norms are equivalent, and $N_1=|||X_1|||$, $\|X_1\|=\|Y_1\|$, 
there exists a constant $C\ge 1$ such that  $C^{-1} N_1\le \|Y_1\|\le C N_1$. \\[0.14cm]
Now let $0<\beta<\alpha$. Hypotheses (i)-(ii) 
show that $\E[(\ln^+ N_1)^\beta]<+\infty$ and $\E[(\ln^- V_1)^\beta]<+\infty$. 
From the previous remarks, we deduce that $\int_S \ell(g)^\beta d\mu(g)<+\infty$ (recall $\mu$ is the law of $Y_1$). 
Denote by $m_\beta$ the previous integral. 
If $\beta\leq1$, then we have $\min\{|t|\ell(g),2\} \leq 2|t|^\beta\ell(g)^\beta$, so that 
$\varepsilon(t)\le 2|t|^\beta m_\beta $. 
If $\beta>1$, then $\varepsilon(t) \leq 2|t|\int\ell(g)d\mu(g)\leq2|t| m_\beta^{\frac{1}{\beta}}$. Thus Theorem III.2 gives 
$\|P_t-P\|_{_{\tiny \cL}}=O(|t|^\beta)$ if $0<\alpha\le 1$, and $\|P_t-P\|_{_{\tiny \cL}}=O(t)$ if 
$1<\alpha\leq2$. Finally, using  $\displaystyle {n\over a_n^\alpha}L(a_n)=1$ and the fact that $L$ is unbounded in the case $\alpha=2$, 
this easily yields the desired statement.  \fdem 

\noindent \noindent{\bf Proof of Theorem I.} Let $(y_n)_{n\geq 1}$ be any sequences of vectors in $\overline{B}$, 
and let $\phi_n$  [resp. $\psi_n$] be the characteristic function of $\frac{\ln\|\widetilde Y^{(n)}y_n\|\, -\, b_n}{a_n}$ 
[resp. of $\frac{\Xi_1+\ldots+\Xi_n\, -\, b_n}{a_n}$]. 
Let $\phi(t) = \mu\otimes\nu (e^{it\xi})$ be the characteristic function of $\xi(Y_1,Z_2)$.
Let $t\in\R$ and $n\in\N^*$ be such that $\frac{t}{a_n}\in I$, and 
set $\ell_n(t) = \Pi(\frac{t}{a_n}){\bf 1}(y_n)$ (cf. Th. III.3). By Theorems III.2-3, one gets $\lim_n\ell_n(t)=1$. 
Furthermore we have 
\begin{eqnarray*}
\phi_n(t) = e^{-i\frac{b_nt}{a_n}}\, \E[e^{i\frac{t}{a_n}\ln\|\widetilde Y^{(n)}y_n\|}] 
&=& e^{-it\frac{b_n}{a_n}}\, P_{\frac{t}{a_n}}{\bf 1}(y_n) \ \ (\mbox{by}\ (\star\star)\ \mbox{in}\ \S\ \mbox{III.1})\\
&=& e^{-it\frac{b_n}{a_n}}\, \lambda(\frac{t}{a_n})^n\, \ell_n(t) +  O(\kappa^n) \ \ (\mbox{by  Th. III.3)}\\
&=& e^{-it\frac{b_n}{a_n}}\, [\phi(\frac{t}{a_n}) + o(\frac{1}{n}))]^n\, \ell_n(t) +  
O(\kappa^n) \ \ (\mbox{Prop. III.1, IV.2}).
\end{eqnarray*}
Finally, since $\psi_n(t) = e^{-it\frac{b_n}{a_n}}\, \phi(\frac{t}{a_n})^n$, one gets
 $\phi_n(t) = \psi_n(t)\, [1 +  o(\frac{1}{n})]^n\, \ell_n(t) +  O(\kappa^n)$, therefore  
$\lim_n\phi_n(t) = \lim_n\psi_n(t)$. Since $(\frac{\Xi_1+\ldots+\Xi_n\, -\, b_n}{a_n})_n$ converges in distribution to a 
stable law of order $\alpha$, the same holds for $(\frac{\ln\|\widetilde Y^{(n)}y_n\|\, -\, b_n}{a_n})_n$. \fdem

%[KLS]  M. A. Kranosel'skij, Je. A. 
 %              Lifshits, A. V. Sobolev (1980). 
  %             { \it Positive Linear Systems - The Method of 
   %            Positive Operators.}
    %           Sigma Series in Applied Mathematics, 
     %          Heldermann Verlag Berlin.

\end{document}